# A fast recurrence for Fibonacci and Lucas numbers

Jeroen van de Graaf

*Departamento da Ciência da Computação*

*Universidade Federal de Minas Gerais*

`jvdg@dcc.ufmg.br`

## Abstract

*We derive the double recurrence $e_n = \frac{1}{2}(a_{n-1} + 5b_{n-1}); f_n = \frac{1}{2}(a_{n-1} + b_{n-1})$ with $e_0 = 2; f_0 = 0$ for the Fibonacci numbers, leading to an extremely simple and fast implementation. Though the recurrence is probably not new, we have not been able to find a reference for it.*

## Introduction

The starting point for this research is the equation

$$\begin{aligned} f_n &= \tfrac{1}{\sqrt{5}}(\tfrac{1+\sqrt{5}}{2})^n - \tfrac{1}{\sqrt{5}}(\tfrac{1-\sqrt{5}}{2})^n \\ &= \tfrac{1}{\sqrt{5}}(1.61803398)^n - \tfrac{1}{\sqrt{5}}(0.61803398)^n \\ &= \tfrac{1}{\sqrt{5}}\varphi^n - \tfrac{1}{\sqrt{5}}\psi^n \end{aligned} \quad (1)$$

Knuth ([GKP89] page 286) pointed out that the second term, $\frac{1}{\sqrt{5}}\psi^n$, is less then $\frac{1}{2}$, so we might just as well compute the first term, and round to the nearest integer. This leads to

$$f_n = \left[ \frac{1}{\sqrt{5}} \left( \frac{1+\sqrt{5}}{2} \right)^n \right] \quad (2)$$

where $[x]$ denotes the rounding function.

The following table shows the results of this approach for small $n$:

$$\begin{array}{c|c|c} n & g_n & f_n \\ \hline 0 & 0.4472 & 0 \\ 1 & 0.7236 & 1 \\ 2 & 1.1708 & 1 \\ 3 & 1.8944 & 2 \\ 4 & 3.0652 & 3 \\ 5 & 4.9597 & 5 \\ 6 & 8.0249 & 8 \\ 7 & 12.9846 & 13 \\ 8 & 21.0095 & 21 \end{array} \quad (3)$$

The problem with this approach is that it does not work well for $n$ large. One needs the exact decimal representation of $\varphi$ and $\psi$ in many decimals, and very precise floating point operations.

## Symbolic algebraical calculations

The main idea of this approach is to not approximate $\sqrt{5}$, but to keep it there as a symbol, and to teach the program the algebraic rules explicitly. For instance, let us define $\alpha = a + b\sqrt{5}$. In analogy with complex numbers, we will call $a$ the real part, and $b$ the algebraic part of $\alpha$. Now

$$\alpha^2 = (a + b\sqrt{5})^2 = (a^2 + 5b^2) + 2ab\sqrt{5} \quad (4)$$

Note that many programming languages today do not impose limit on the size of integer arithmetic, so we will not lose precision.

As second idea results from computing $f_2$ using the first fórmula:

$$
\begin{aligned}
f_2 &= \tfrac{1}{\sqrt{5}}(\tfrac{1+\sqrt{5}}{2})^2 - \tfrac{1}{\sqrt{5}}(\tfrac{1-\sqrt{5}}{2})^2 \\
&= \tfrac{1}{4\sqrt{5}}(\cancel{(1+5)} + 2\sqrt{5} - (\cancel{(1+5)} - 2\sqrt{5})) \\
&= \tfrac{1}{4\sqrt{5}}(4\sqrt{5})) \\
&= 1
\end{aligned}
\tag{5}
$$

Note how, in the denominators, the real parts cancel. This must be so because if not, the factor $\tfrac{1}{\sqrt{5}}$ would cause $f_n$ to have an algebraic part not equal to 0, whereas we know that $f_n$ is an integer. Note also how the algebraic parts of the denominators double, yielding $f_n$ after division by (or multiplication with) some rational multiple of $\sqrt{5}$.

From these calculations we deduce that if

$$
\psi^n = \left(\frac{1-\sqrt{5}}{2}\right)^n = a - b\sqrt{5} = \alpha
\tag{6}
$$

then

$$
\varphi^n = \left(\frac{1+\sqrt{5}}{2}\right)^n = a + b\sqrt{5} = \overline{\alpha}
\tag{7}
$$

and that

$$
\begin{aligned}
f_n &= \tfrac{1}{\sqrt{5}}\varphi^n - \tfrac{1}{\sqrt{5}}\psi^n \\
&= \tfrac{1}{\sqrt{5}}(a + b\sqrt{5}) - \tfrac{1}{\sqrt{5}}(a - b\sqrt{5}) \\
&= \tfrac{1}{\sqrt{5}}(a - a) - \tfrac{1}{\sqrt{5}}((b + b)\sqrt{5}) \\
&= 2b
\end{aligned}
\tag{8}
$$

In other words, by computing $\psi^n$ we get $\varphi^n$ almost for free, and we can compute $f_n$ right away.

However, a problem with

$$
f_n = \tfrac{1}{\sqrt{5}}(\tfrac{1+\sqrt{5}}{2})^n - \tfrac{1}{\sqrt{5}}(\tfrac{1-\sqrt{5}}{2})^n
\tag{9}
$$

are the factors $\tfrac{1}{2}$ in $\varphi$ e $\psi$. One possibility is to compute

$$
f_n = \tfrac{1}{\sqrt{5}} \cdot \tfrac{1}{2^n}((1+\sqrt{5})^n - (1-\sqrt{5}))^n
\tag{10}
$$

but then the powers of $(1 \pm \sqrt{5})$ will grow very fast. It is more elegant to change our algebraic representation, incorporating $\tfrac{1}{2}$.

So we now change our representation. Instead of associating the pair $\alpha = (a, b)$ to $a + b\sqrt{5}$, we define

$$
(a, b) \longleftrightarrow \tfrac{1}{2}(a + b\sqrt{5})
\tag{11}
$$

Obviously $\varphi \longleftrightarrow (1, 1)$ and $\psi \longleftrightarrow (1, -1)$, a small caveat is that $1 \longleftrightarrow (2, 0)$. In this new notation, we obtain for instance

$$
\alpha^2 = (\tfrac{1}{2}(a+b\sqrt{5}))^2 = \tfrac{1}{4}((a^2 + 5b^2) + 2ab\sqrt{5}) = \tfrac{1}{2}(\tfrac{1}{2}(a^2 + 5b^2) + ab\sqrt{5}) = (a\prime, b\prime) = (12(a2 + 5b2), ab5)
\tag{12}
$$

so $(a', b') = (\tfrac{1}{2}(a^2 + 5b^2), ab\sqrt{5})$. Note that the division by 2 of the real component results in an integer as long as $a$ e $b$ have the same parity, which will always be the case.

In order to compute $\psi^n = \psi^{n-1} \cdot \psi$, we must multiply an arbitrary number $\alpha = \tfrac{1}{2}(a + b\sqrt{5})$ by $\psi = \tfrac{1}{2}(1 - \sqrt{5})$ in the new representation. So

$$
\tfrac{1}{2}(a + b\sqrt{5}) \cdot \tfrac{1}{2}(1 - \sqrt{5}) = \tfrac{1}{2}\left(\tfrac{1}{2}(a - 5b) + \tfrac{1}{2}(b - a)\sqrt{5}\right) = \tfrac{1}{2}(a' + b'\sqrt{5})
\tag{13}
$$

so

$$
\begin{cases} a' = \tfrac{1}{2}(a - 5b) \\ b' = \tfrac{1}{2}(b - a) \end{cases}
\tag{14}
$$

From this we can derive a recurrence, recalling that $\psi = \tfrac{1}{2}(1 - \sqrt{5}) \longleftrightarrow (a_1 = 1, b_1 = -1)$:

$$
\begin{cases} a_1 = 1;\ b_1 = -1 \\ a_n = \tfrac{1}{2}(a_{n-1} - 5b_{n-1}) \end{cases}
\tag{15}
$$

$$\left\{ b_n = \tfrac{1}{2}(b_{n-1} - a_{n-1}) \right.$$

where $f_n = -b_n$.

After some reflection on the minus-signs in this formula, it becomes clear that they can all be changed by plus-signs. We can compute $\varphi^n$ instead of $\psi^n$, we were too cautious.

## The main recurrence

So this leads us to the following recurrence, which is the main result of this note:

$$\begin{cases} e_0 = 2; \quad f_0 = 0 \\ e_n = \tfrac{1}{2}(e_{n-1} + 5f_{n-1}) \\ f_n = \tfrac{1}{2}(e_{n-1} + f_{n-1}) \end{cases} \quad (16)$$

Though surely some mathematician must have found this result already, we did not find a reference to this formula in the literature or on the internet. (If you know a reference, please contact the author.) Interestingly, the successive values $e_i$ also follow a Fibonacci recurrence but with different initial values, resulting in a sequence known as the Lucas numbers: $2, 1, 3, 4, 7, 11, 18, 29, \ldots$

In Python3, the following code prints the first 100 Fibonacci numbers.

```
(e,f)=(2,0)
for i in range(100):
    print(i,f)
    e,f = (e+5*f)//2, (e+f)//2
```

## Square-and-Muliply

Our recurrence is clearly more complicated then the original formula $f_n = f_{n-1} + f_{n-2}$, but, unlike latter, our new recurrence allows us to skip intermediate Fibonacci numbers, which allows us to implement very efficient algorithms for $n$ large.

As a starter, note that $\varphi^{64} = (((((\varphi^2)^2)^2)^2)^2)^2$. So if we teach the program how to square a pair $(a,b) \longleftrightarrow \tfrac{1}{2}(a + b\sqrt{5})$, we can compute $f_{2^m}$ in $m$ steps. One easily verifies that

$$(\tfrac{1}{2}(a + b\sqrt{5}))^2 = (a^2 + 5b^2)/2 + ab\sqrt{5} \quad (17)$$

For instance, the following Python3 program computes some $f_n$ for $n$ a power of 2.

```
a,b=(1,1)
for m in range(1,8):
    a,b=(a*a+5*b*b)//2,(2*a*b)//2
    print(2**m,b)
```

Generalizing beyond powers of $2$, a well-known algorithm for optimizing exponentiations is *Square-and-Multiply*, which computes $a^n$ in $\lg n$ iterations, taking advantage of the binary representation of the exponent $n$. See for instance section 14.6.1 in [MvOV97]. For example, $49 = 110001_2$, and

$$a^{49} = a^{48} \cdot a = ((a^3)^{16}) \cdot a = (((((a^2 \cdot a)^2)^2)^2)^2) \cdot a \quad (18)$$

Roughly speaking, *Square-and-Multiply* uses some variable $x$ which starts as $a$, which squares in each iteration, and multiplies $x$ by $a$ if there is a `1` in the binary representation of $n$.

## Main algorithm for Fibonacci

Straightforward application of *Square-and-Multiply* to our problem yields the following program.

```
1  def fib(n):
2      bits_n = bin(n)[2:] # get binary representation
3      a,b=(2,0)
4      for i in bits_n:
5          a,b=(a*a+5*b*b)//2,(2*a*b)//2          # "square"
6          if i == '1': a,b = (a+5*b)//2, (a+b)//2 # "multiply"
7      return b
```

This program runs very fast; it determines $f_{1000000}$ in less than 1 second in the Jupyter environment on a off-the-shelf laptop.

This new method takes $\lg n$ squarings and $w(n)$ multiplications, where $w(n)$ is the Hamming weight of $n$, i.e. the number of 1s in the representation of $n$. However, we are dealing with operations on pairs $(a,b)$ which use aditions and multiplications on integers which expand in size. Observe that a multiplication of the pair $(a,b)$ by $\varphi$

$$\begin{cases} a' = \frac{1}{2}(a+5b) \\ b' = \frac{1}{2}(b+a) \end{cases} \quad (19)$$

can be implemented with additions, subtractions and shifts only (since $5b = 4b+b$, and multiplying by 4 is a shift to the left of 2 positions).

Squaring a pair $(a,b)$ involves 2 squaring and 1 multiplication of integers.

$$\begin{cases} a' = a^2 + 5b^2 \\ b' = ab \end{cases} \quad (20)$$

## Comparison with other programs

We have not had the time to make detailed performance analysis, but compared with the algoritms presented in [Dasdan18] our algorithm is much simpler. It also seems much faster.

It seems that in that paper there are two approaches which compete. The first approach is based on recurrences of the type

$$\begin{aligned} f_{2n-1} &= f_n{}^2 + f_{n-1}{}^2 \\ f_{2n} &= (2f_{n-1} + f_n)f_n. \end{aligned} \quad (21)$$

Though this also leads to $\lg n$ iterations and has 2 resp. 1 multiplication, note that the recursion is more complicated: for each $f_{2n}$ or $f_{2n-1}$, we need *two* values, $f_{n-1}$ and $f_n$ . Compare this to the right-to-left recurrence of *Square-and-Multiply*,

$$\begin{aligned} \varphi_{2n+1} &= (\varphi_n)^2 \cdot \varphi \\ \varphi_{2n} &= (\varphi_n)^2 \end{aligned} \quad (22)$$

in which only one value is needed. Moreover, these two recursions tend to spread out, so as to avoid wasteful recursion calls, one needs to implement memoization (once a Fibonacci number has been computed, it is stored in 'cache' to avoid that the same calculation takes place later). To put it differently: yes, the complexity of both is $O(\lg n)$, but the constant in our algoritm is much better.

In the second approach the Fibonacci numbers appear as powers of the matrix $M = \begin{pmatrix} 1 & 1 \\ 1 & 0 \end{pmatrix}$ since

$$M^n = \begin{pmatrix} 1 & 1 \\ 1 & 0 \end{pmatrix}^n = \begin{pmatrix} f_n & f_{n-1} \\ f_{n-1} & f_{n-2} \end{pmatrix}^n \quad (23)$$

If we look at the operations for *Square-and-Multiply* for $M$, we see that the multipliction of an arbirary matrix $\begin{pmatrix} a & b \\ b & c \end{pmatrix}$ by $M$ results in

$$\begin{pmatrix} a & b \\ b & c \end{pmatrix} \cdot \begin{pmatrix} 1 & 1 \\ 1 & 0 \end{pmatrix} = \begin{pmatrix} a+b & a \\ b+c & b \end{pmatrix} = \begin{pmatrix} a+b & a \\ a & b \end{pmatrix} \quad (24)$$

since in this case $a = b + c$.

The squaring of an arbitrary matrix $\begin{pmatrix} a & b \\ b & c \end{pmatrix}$ results in

$$\begin{pmatrix} a & b \\ b & c \end{pmatrix} \cdot \begin{pmatrix} a & b \\ b & c \end{pmatrix} = \begin{pmatrix} a^2 + b^2 & ab + bc \\ ab + bc & b^2 + c^2 \end{pmatrix} \tag{25}$$

corresponding to 3 integer squarings and 2 multiplications per operation. Earlier we reported 2 squarings and 1 multiplication, so the new algorithm seems to be faster. However, a more detailed comparison would be necessary.

## Concluding notes

- This draft is a preliminary version, and likely to contain inperfections. For instance, it was typeset quick-and-dirty with Typora. I will convert to Latex ASAP.
- Alejandro Hevia suggest the link to [OEIS]. I still need to look at the references, there are **many**.
- Why people want to know $f_{1000000}$? It beats me, but it can be done :-).

## References

[GKP] Graham, Knuth, Parashnik. *Concrete mathematics*

[D2018] Dasdan, A. *Twelve simple algoritms to compute Fibonacci numbers* https://arxiv.org/pdf/1803.07199.pdf

[MvOV97] Menezes, van Oorschot, Vanstone *Handbook of Applied Cryptography*

https://cacr.uwaterloo.ca/hac/

[OEIS] https://oeis.org/ search for "fibonacci".